\documentclass[a4paper, 12pt]{article}
\usepackage{amsmath}
\usepackage{amsfonts}
\usepackage{amssymb}
\usepackage{amscd}
\usepackage{amsthm}
\usepackage[all,cmtip]{xy}

\newfont{\cyrr}{wncyr10}
\newcommand\Sh{\mbox{\cyrr Sh}}

\newcommand{\da}{\mathbb{A}}
\newcommand{\Z}{{\mathbb Z}}
\newcommand{\Fq}{{\mathbb F}_q}
\newcommand{\Q}{{\mathbb Q}}

\newcommand{\R}{{\mathbb R}}
\newcommand{\C}{{\mathbb C}}

\newcommand{\Fo}{{\rm F}}
\newcommand{\Co}{{\rm C}}
\newcommand{\Fb}{{\mathbb F}}

\newcommand{\T}{{\mathbb T}}

\newcommand{\A}{{\mathbf A}}

\newcommand{\res}{{\rm res}}
\newcommand{\Aut}{{\rm Aut}}

\newcommand{\Hom}{{\rm Hom}}

\newcommand{\Gal}{{\rm Gal}}

\newcommand{\Pic}{{\rm Pic}}

\newcommand{\Ker}{{\rm Ker}}

\newcommand{\rk}{{\rm rk}}

\newcommand{\Tr}{{\rm Tr}}

\newcommand{\Spec}{{\rm Spec}}

\newcommand{\Det}{{\rm Det}}

\renewcommand{\dim}{{\rm dim}}

\newcommand{\GL}{{\rm GL}}

\newcommand{\Div}{{\rm Div}}

\theoremstyle{plain}

\newtheorem{defin-prop}[theor]{Definition-Proposition}

\theoremstyle{remark}

\theoremstyle{definition}

\title{Representations  of Higher Adelic Groups and Arithmetic
\footnote{To appear in Proceedings of the International Congress of Mathematicians, Hyderabad, India, 2010, vol. 1.}
}
\author{A. N. Parshin\footnote{Steklov Mathematical Institute, Russian Academy of Sciences, Gubkina str 8, 119991 Moscow, Russia. I am very thankful to Lawrence Breen,   Denis V. Osipov,  Vladimir L. Popov  and  Yuri G. Zarhin  who have read the text and made many valuable remarks. The author was supported by RFBR (grant no.
08-01-00095-a), by a program  for supporting
Leading Scientific Schools (grant no. Nsh-4713.2010.1
).}}
\date{}

\begin{document}

\maketitle

\bigskip

What do we mean by local ? To get an answer to this  question let us start from  the following two problems.

First problem is from number theory. When does the diophantine equation
$$
f(x, y, z) = x^2 - ay^2 -bz^2 = 0, \quad a, b, ~\in \Q^*
$$
have a non-trivial solution in rational numbers ? In order to solve the problem, let us  consider the quadratic norm residue symbol
$(-,-)_p$ where $p$ runs through all primes $p$ and also $\infty$. This symbol is a bi-multiplicative map $(-,-)_p : \Q^* \times \Q^* \rightarrow
\{\pm 1\}$ and it is easily computed in terms of the Legendre symbol. Then, a non-trivial solution exists if and only if,  for any $p$,  $ (a, b)_p = 1$. However, these conditions are not independent:
\begin{equation}
\label{1}
\prod_p (a, b)_p = 1.
\end{equation}
This is essentially the Gauss  reciprocity law in the Hilbert form.

The ``points" $p$ correspond to {\em all possible} completions of the field $\Q$ of rational numbers, namely to
the $p$-adic fields $\Q_p$ and the field $\R$ of real numbers. One can show that the equation $f = 0$ has a non-trivial solution
in $\Q_p$ if and only if $(a, b)_p = 1$.

The second  problem comes  from complex analysis. Let $X$ be a compact Riemann surface (= complete smooth algebraic curve defined over $\C$). For a point  $P \in X$, denote by $K_P = \C((t_P))$ the field of  Laurent formal power series in a local coordinate $t_P$  at the point $P$. The field  $K_P$
contains the ring $\widehat{\cal O}_P = \C[[t_P]]$ of  Taylor formal power series. These have an invariant meaning and are  called the local field and the local ring  at  $P$ respectively.
 Let us now fix finitely many points $P_1,\dots , P_n \in X$ and assign to every $P$ in $X$  some   elements $f_P$ such that $f_{P_1} \in K_{P_1},\dots , f_{P_n}\in K_{P_n}$ and $f_P = 0$ for all other points.

When does there exist a meromorphic (=rational) function $f$ on $X$ such that
\begin{equation}
\label{2}
 f_P - f \in \widehat{\cal O}_P \quad \mbox{for every}~P \in X?
\end{equation}
The classical answer to this Cousin problem is the following: there exists such an $f$ whenever for any regular differential form $\omega$ on $X$
\begin{equation}
\label{3}
\sum_P \res_P(f_P\omega) = 0.
\end{equation}
The space of regular differential forms has dimension $g$ (= genus of $X$) and in this way one gets finitely many conditions on the data $(f_P)$. The residue is an additive map $\res_P : \Omega^1(K_P) \rightarrow \C$ and is easily computed in terms of the local decomposition of the differential form $\omega \in   \Omega^1(K_P)$.     Note that ``locally",  problem (\ref{2}) can be solved
for any point $P$. Behind our global conditions (\ref{3}), we have the following residue relation:
\begin{equation}
\label{4}
\sum_P \res_P(\eta) = 0
\end{equation}
for any {\em meromorphic} differential form $\eta$  on $X$.

We see some similarity between these two problems, which  belong to very different parts of our science. The explanation lies in the existence of a very
deep analogy between numbers and functions, between number fields and fields of algebraic functions. This analogy goes back to the nineteenth century,
possibly to Kronecker. The leading role in  the subsequent development belongs to Hilbert. The  analogy was one of his beloved ideas, and  thanks to Hilbert  it  became one of the central ideas in the development of number theory during the twentieth  century. Following this analogy, we can compare algebraic curves over $\C$ (= compact Riemann surfaces) and number fields (= finite extensions of $\Q$). In particular, this includes a comparison of  {\em local} fields such as that between the fields $\C((t))$ and $\Q_p$. Their similarity was already pointed out by Newton\footnote{He compared the power series and the expansions of rational numbers in powers of $p$~ (for $p = 10$).}.

In modern terms, we have two kinds of geometric objects. First, a complete algebraic curve $X$, containing an affine curve $U = \Spec(R)$\footnote{Here, $\Spec(R)$ is the set of prime ideals in the ring $R$ together with the additional structure of a scheme.}, where $R$ is the ring of regular functions on $U$:
$$
\mbox{(geometric picture)}\quad X \supset U\quad \mbox{and  finitely many points}\quad P \in X.
$$
Next, if we turn to  arithmetic,  we have a finite extension $K \supset \Q$ and the ring $R \subset K$ of integers. We write
$$
\mbox{(arithmetic picture)}\quad X \supset U = \Spec(R) \quad \mbox{and finitely many infinite places}~ P \in X.
$$
The places (``points") correspond to the embeddings of $K$ into the fields $\R$ or $\C$. Here, $X$ stands for the as yet  not clearly defined complete ``arithmetical''
curve, an analogue of the curve $X$ in the geometric situation. The analogy between both $U$'s is very clear and transparent. The rings $R$ are the Dedekind rings of the Krull dimension\footnote{That is, the length of a maximal chain of prime ideals. The ring $R$ itself is not a prime ideal.}  1.
The nature of the additional points (outside $U$) are  more complicated. In the geometric case, they also  correspond to the non-archimedean valuations on the curve $X$, whereas in the arithmetical case these  infinite places are a substitute for the archimedean valuations of the field $K$.

In algebraic geometry, we also have  the theory of algebraic curves defined over a finite field $\Fb_q$ and this theory,  being arithmetic in its nature, is much  closer to the theory of number fields than the theory of algebraic curves over $\C$. The main construction on both sides of the analogy is the notion of a  local field. These local fields appear  into the following table:
\vskip 0.5cm
\begin{center}
\begin{tabular}{|c|c|c|}
\hline
 $\mbox{dimension}$ & geometric case & arithmetic case \\
 $> 2$     &          ...            &                 ...        \\
\hline
2  & ? & $\R((t)),~\C((t))$ \\
\hline
1    &  $\Fb_q((t))$      &  $\Q_p,~\R,~\C$ \\
\hline
 0    &        $\Fb_q$     &   $\Fb_1$ \\
\hline
\end{tabular}
\end{center}
\vskip 0.5cm
Here $\Fb_1$ is the so-called  ``field'' with  one element, which is quite popular nowadays. We will see soon   why the fields $\R((t))$ and $\C((t))$ belong to the higher level of the table than the fields $\Q_p$  or $\R$. More on the analogy between geometry and arithmetic can be found in \cite{P9}.

\section{$n$-dimensional Local Fields and Adelic Groups}

Let us consider  algebraic varieties $X$ (or Grothendieck schemes) of dimension greater than one. It appears that we have a
well established notion of something local attached to a point $P \in X$. One can take a neighborhood of $P$,  e.g. affine, complex-analytic if
$X$ is defined over $\C$, formal and so on. In this talk we will advocate the  viewpoint that the genuine local objects on the varieties are not the points with some    neighborhoods but the maximal ordered  sequences (or flags) of  subvarieties, ordered by inclusion.

 If $X$ is a variety  (or a scheme) of dimension $n$
and
$$X_0\subset X_{1}\subset\dots X_{n-1}\subset X_n=X$$
is a flag of irreducible subvarieties ($\dim (X_i)=i$)
then one can define a certain ring
$$K_{X_0,\dots,X_{n-1}}$$ associated to the flag. In the case where
  all the subvarieties are
regularly embedded, this ring is an $n$-dimensional local field.

\smallskip

{\bf Definition 1}.
 Let $K$ and $k$ be fields. We say that $K$ has a structure
 of an {\it $n$-dimensional local field} with the {\it  last residue field} $k$ if either
 $n=0$ and $K=k$ or $n\geq 1$ and $K$ is the fraction field of a complete discrete
 valuation ring ${\cal O}_K$ whose residue  field $\bar K$ is a local field
 of dimension $n-1$ with the last residue field $k$.

\smallskip

Thus, an $n$-dimensional local field
  has the following inductive structure:
$$ K=:K^{(0)}\supset {\cal O}_K\to {\bar K}=:K^{(1)}\supset
   {\cal O}_{\bar K}\to {\bar K}^{(1)}=:K^{(2)}
   \supset {\cal O}_{K^{(2)}}\to\ldots\to {\bar K}^{(n)}=k            $$
where ${\cal O}_F$ denotes the valuation ring of the valuation on
$F$ and $\bar{F}$ denotes the residue field.

The  simplest example of an  $n$-dimensional local field is the field
$$
K = k((t_1))((t_2))\dots((t_n))
$$
of iterated Laurent formal power series.  In dimension one, there are examples from the
 table.  However, fields such as  $\R$ or $\C$ are not  covered by this definition. Concerning classification of the local fields
 see \cite{FP}.

One can then form the adelic group (actually, the  ring)
$$\Bbb A_X= {\prod}' K_{X_0,\dots,X_{n-1}}$$
where the product is taken
over all the flags  with respect to certain restrictions on components
of adeles. For  schemes over a finite field $\Fq$, this  is
the ultimate definition of the adelic space attached to $X$. In
general, one must  extend it   by adding archimedean components, such as  the fields
$\R$ or $\C$ in dimension one.

In dimension one, the local fields and  adelic groups
are   well-known
tools of  arithmetic. They were introduced by C. Chevalley in the
1930s and were used to formulate and   solve
many  problems in number theory and algebraic geometry
(see, for example, \cite{A,W2}). These constructions are associated
with  fields of algebraic numbers and  fields of algebraic
functions in  one variable over a finite field, that is  with schemes of dimension $1$. A need for
such constructions  in higher dimensions was realized by the author  in the  1970s.
They were  developed
in  the local case for any dimension and  in the global case for
dimension two \cite{P1, P2, FP, P6}. This approach was extended  by A. A. Beilinson
to the schemes of an arbitrary dimension \cite{Be, Hu}.
In this talk, we restrict ourselves to the case of dimension two.

Let $X$ be a smooth  irreducible surface over a field $k$ (or an arithmetic surface), let
 $P$ be a closed point of $X$ and let  $C\subset X$ be an irreducible
curve such that $P\in C$. We denote by ${\cal O}_{X,P}$ the local ring at the point $P$, that is the ring of rational functions which are
regular at $P$. We denote also by ${\cal O}_C$ the  ring of rational functions on $X$ which have no pole along the $C$.

If  $X$ and $C$ are smooth at $P$, then we pick  a local
equation $t\in {\cal O}_{X,P}$ of $C$ at $P$  and  choose
$u\in {\cal O}_{X,P}$  such that $u|_C\in {\cal O}_{C,P}$
is a local parameter at $P$. Denote by $\wp$ the ideal  in ${\cal O}_{X,P}$ defining the
curve $C$ near $P$. We can introduce a two-dimensional local field
$K_{P,C}$ attached
to the pair $P, C$ by the following procedure  which includes completions and
localizations:
$$
\begin{matrix}
{\cal O}_{X,P}&&\\
\vert &  & \\
{\widehat{\cal O}}_{X,P} & = & k(P)[[u,t]] \supset \wp = (t)\\
\vert &  & \\
({\widehat{\cal O}}_{X,P})_{\wp} & = &
\text{discrete valuation ring with residue field $k(P)((u))$}\\
\vert &  & \\
{\widehat{\cal O}}_{P,C}:= \widehat{({\widehat{\cal O}}_{X,P})}_{\wp} & = & k(P)((u))[[t]]\\
\vert &  & \\
K_{P,C}:=\text{Frac}\,({\widehat{\cal O}}_{P,C}) & =& k(P)((u))((t))
\end{matrix}
$$
Note that the left-hand  construction is meaningful even
{\it without}  smoothness of the curve $C$ (it is sufficient to assume that $C$ has only one formal branch near $P$). In the general  case, the ring  $K_{P,C}$ is a finite direct sum of  2-dimensional local fields. If  $P$ is smooth then the field $K_{P,C}$ has the following informal interpretation. Take a function $f$ on $X$. We can, first, develop $f$ as a formal power series
in the variable $t$ along the curve $C$ and then every coefficient of the series restricted to $C$ can be further developed
as a formal power series in the variable $u$. The local  field $K_{P,C}$ is  a kind of completion of the field of rational functions $K = k(X)$ on $X$.
It carries a discrete valuation $\nu_C : K_{P,C}^* \rightarrow \Z$ defined by the powers of the ideal $\wp$.

Let  $K_P$ be the minimal subring of $K_{P,C}$
which contains both $k(X)$ and ${\widehat{\cal O}}_{X, P}$.
In general, the ring $K_P$ is not a field.
Then $K\subset K_P\subset K_{P,C}$ and  there is another intermediate
subring $K_C =\text{Frac}\,(\widehat{\cal O}_C) \subset K_{P,C}$.
 We can compare the structure
of the local adelic components in dimensions one and two:
$$
\xymatrix@!0{
K_P\ar@{-}[dd] & & &  & K_{P,C}\ar@{-}[dl] \ar@{-}[dr]  &
\\
 & & & K_P\ar@{-}[dr]  & & K_C \ar@{-}[dl]
\\
K & & & & K &
}
$$

The global adelic group is a certain subgroup of the ordinary  product of all  two-dimensional local fields. Namely, a collection  $( f_{P, C})$  where
$f_{P, C} \in K_{P, C}$ belongs to  $\Bbb A_X$ if the  following two conditions are satisfied:
\begin{itemize}
\item
$$ \{ f_{P, C}\} \in \da_C((t_C))$$
for a fixed irreducible curve $C \subset X$ and  a local equation $t_C=0$ of the curve $C$ on some open affine subset $U \subset X$ and
\item we have $\nu_C ( f_{P, C}) \ge 0$, or equivalently
$$ \{ f_{P, C}\} \in \da_C[[t_C]],
$$
 for  all but   finitely many  irreducible curves $C \subset X$.
\end{itemize}
Here we reduced the definition of the adelic group to the classical case of algebraic curves $C$. Recall that a collection  $(f_P,~P \in C)$ belongs to the adelic (or restricted) product $\da_C$ of the local fields $K_P$ if and only if for almost all points $P$ we have $f_P \in \widehat{\cal O}_P$.

\smallskip

What can one  do with this  notion of the local field and why is it   really local ? To get some understanding of this,  we would like to develop the
above examples (of residues and symbols) in  dimension two. For any flag $P\in C$ on a surface $X$ and a rational differential form $\omega$ of degree 2
we can define the residue
$$
\res_{P, C}(\omega) = \Tr_{k(P)/k} (a_{-1, -1})
$$
where $\omega = \sum_{i,j}a_{i, j}u^it^jdu\wedge dt$ in the field $K_{P,C}\cong k(P)((u))((t))$.
Then, instead of the simple relation (4) on an algebraic curve,  we get
two types of relations on the projective surface $X$ \cite{P2}
\begin{equation}\label{res1}
 \sum_{P \in C} \res_{P, C}(\omega) =  0, \quad \mbox{for any fixed  curve}~C,
\end{equation}
\begin{equation}\label{res2}
 \sum_{C \ni  P} \res_{P, C}(\omega) = 0, \quad \mbox{for any fixed point}~P.
\end{equation}

At the same time,  we can define certain  symbols (bi-multiplicative and three-multiplicative) \cite{P1}
$$
(-, -)_{P,C}: K_{P,C}^*\times K_{P,C}^* \rightarrow \Z \quad \mbox{and}~ (-, -, -)_{P,C}: K_{P,C}^*\times K_{P,C}^* \times K_{P,C}^* \rightarrow k^*$$
which are   respectively generalizations of the valuation $\nu_P : K^*_P \rightarrow \Z$ and the norm residue symbol $(-, -)_P : K_{P}^*\times K_{P}^*$ (actually, the tame symbol)  on an  algebraic curve $C$.
The reciprocity laws have the same structure as the residue relations. In particular, if $f, g, h \in K^*$ then
$$
 \prod_{P \in C} (f, g, h)_{P,C} = 1, \quad \mbox{for any fixed curve}~C,
$$
$$
 \prod_{C \ni  P} (f, g, h)_{P,C} = 1, \quad \mbox{for any fixed point}~P.
$$

This shows that in dimension two there is a symmetry between points $P$
and curves $C$ (which looks like  the classical duality between points and
lines in projective geometry).

If $C$ is a curve then the space $\da_C$ contains the important subspaces
$\da_0 = K = k(C)$ of principal adeles (rational functions diagonally embedded into the adelic group) and $ \da_1 = \prod_{P \in C}
{\widehat{\cal O}}_P$ of integral adeles.  These  give rise to   the adelic
complex
\begin{equation}\label{ad}
 \da_0 \oplus \da_1 \rightarrow \da_C .
\end{equation}
This complex computes the cohomology of the structure sheaf ${\cal O}_C$. If $D$ is a divisor on $C$ then the cohomology of the sheaf
${\cal O}_C(D)$ can be computed using the adelic  complex (\ref{ad}) where the subgroup  $\da_1$ is replaced by the subgroup $\da_1(D) =  \{ (f_P) \in A_C :
\nu_P (f_P) + \nu_P(D) > 0~\mbox{for any}~P \in C$\}.

In dimension two,
there is a much more complicated structure of subspaces in $\da_X$ (see \cite{P6}).
Among the others, it includes three subspaces $\da_{12} = {\prod}'_{P\in C}
\widehat{\cal O}_{P,C}$, $\da_{01} = {\prod}'_{C \subset X}
K_C$  and $\da_{02} = {\prod}'_{P\in X}
K_P$.  We set $ \da_0 = \da_{01}\cap \da_{02},  \da_1 = \da_{01}\cap \da_{12}$ and $ \da_2 = \da_{02}\cap \da_{12}$,
and  arrive at  an adelic complex
$$
 \da_0 \oplus    \da_1 \oplus  \da_2 \rightarrow  \da_{01}\oplus \da_{02} \oplus  \da_{12} \rightarrow  \da_X.
$$
Once again, the complex computes the cohomology of the sheaf  ${\cal O}_X$.
One can  extend these complexes to the case of arbitrary schemes $X$ and any coherent sheaf on $X$
(see \cite{Be, Hu, FP}).

The last issue which we will discuss in  this section  is the relation between  the residues and  Serre duality for coherent sheaves. We will only consider the construction of the fundamental class  for the sheaf of differential forms. For curves $C$, we have an isomorphism $H^1(C, \Omega^1_C) \cong \Omega^1(\da_C)/\Omega^1(\da_0) \oplus \Omega^1(\da_1)$. The fundamental class isomorphism  $H^1(C, \Omega^1_C) \cong k$ can be defined as the sum of residues on  $\Omega^1(\da_C)$. The residues relation (3) shows that this sum vanishes on the subspace $\Omega^1(\da_0)$ (and it vanishes on the other subspace $\Omega^1(\da_1)$  for trivial reasons). The same reasoning works in the case of surfaces. We have an isomorphism
$$
H^2(C, \Omega^2_X) \cong \Omega^2(\da_X)/\Omega^2(\da_{01}) \oplus \Omega^2(\da_{02})\oplus \Omega^2(\da_{12}) \rightarrow k,
$$
where the last arrow is again the sum of residues over all flags $P \in C \subset X $. The correctness of this definition follows from the residues relations (\ref{res1}) and (\ref{res2}). We refer to \cite{P2, Be, FP} for the full description of the duality.

\section{Harmonic Analysis on Two-dimensional Schemes}
In  the 1-dimensional  case,  local fields and  adelic groups both carry a natural topology for which they are locally compact groups and
classical harmonic analysis on  locally compact groups can  therefore be applied to  this situation. The study of representations of algebraic groups over local fields and adelic groups is a broad subfield of  representation theory, algebraic geometry and number theory. Even for abelian groups, this line of thought has very
nontrivial applications in  number theory, particularly to the study of L-functions of  one-dimensional schemes (see below).
The first preliminary step is the existence of a  Haar measure on  locally compact groups.
The  analysis starts with a definition of certain  function spaces.

We have two sorts  of  locally compact groups. The groups of the first type are totally disconnected such as the fields $\Q_p$ or $\Fq ((t))$. These groups are related with varieties defined over a finite field. The groups of the second type are  connected Lie groups such as the fields $\R$ or $\C$.

If  $V$ is a locally compact abelian group of the first type let us consider the following spaces of functions (or distributions) on $V$:
$$
\begin{array}{rcl}
{\cal D}(V) &=& \{\mbox{locally constant functions with compact support}\}\\
\tilde{{\cal E}}(V) &=& \{\mbox{uniformly locally constant functions}\}\\
{\cal E}(V) &=& \{\mbox{all locally constant functions}\}\\
{\cal D}'(V) &=& \{\mbox{the dual to}~{\cal D}(V) \mbox{, i.e. all distributions}\}\\
\tilde{{\cal E}}'(V) &=& \{\mbox{the  ``continuous'' dual
to}~{\tilde{\cal E}(V)}\}\\
{\cal E}'(V)&=& \{\mbox{the ``continuous'' dual to} ~{\cal E}(V) \mbox{,
i.e. distributions with compact support}\} \mbox{.}
\end{array}
$$
These are the classical spaces introduced by F. Bruhat \cite{Br} and the more powerful  way to develop the harmonic
analysis is the categorical point of view.
First, we need  definitions of direct and inverse images  with respect to the continuous  homomorphisms.

Let  $f: V \rightarrow W$ be a strict homomorphism\footnote{This means that $f$ is a composition of an open epimorphism and a closed monomorphism.} of locally  compact groups $V$ and $W$.
Then the inverse image     $f^*: {\cal D}(W) \rightarrow {\cal D}(V)$ is defined  if and only if the kernel of $f$ is compact.
The direct image   $f_*: {\cal D}(V) \otimes \mu (V) \rightarrow {\cal D}(W)$ is defined if and only if the cokernel of $f$ is discrete.
Here, $\mu (V)$ is a (1-dimensional) space
of  Haar measures on $V$. For the spaces like ${\cal E}, \tilde{\cal E}$ the inverse image is defined for any $f$, but the direct image is defined
if and only if   the kernel is compact and the cokernel  is discrete. For the distribution spaces the corresponding conditions are the dual ones.
Therefore, we see  that these
maps do not exist for arbitrary homomorphisms in our category and  there are  some ``selection rules''.

The Fourier transform  $\Fo$ is defined as
a map from ${\cal D}(V) \otimes \mu (V)$ to
${\cal D}(\check V)$ as well as for the other
types of spaces. Here,  $\check V$
is the dual group. The main result is the
following Poisson formula
$$
\Fo(\delta_{W,\mu_0} \otimes \mu) =
\delta_{W^{\perp},\mu^{-1}/\mu^{-1}_0}
$$
for any closed subgroup $i:W \rightarrow V$.
Here $\mu_0 \in \mu(W)\subset {\cal D}'(W), \mu \in
\mu(V)\subset {\cal D}'(V), \delta_{W,\mu_0} = i_*(1_W \otimes
\mu_0)$ and $W^{\perp}$ is the annihilator of
$W$ in $\check V$.

This general formula is very efficient when applied to
the {\em self-dual} (!)  group $\da_C$.
 The standard subgroups in $\da_C$ have their characteristic functions
$ \delta_{\da_1(D)} \in  {\cal D}(\da_C)$ and $ \delta_K \in  {\cal D}'(\da_C)$ .
We have
\begin{equation}\label{p1}
\Fo(\delta_{\da_1(D)}) = \mbox{vol}(\da_1(D))\delta_{\da_1((\omega)-D)},
\end{equation}
\begin{equation}\label{p2}
 \Fo(\delta_K) = \mbox{vol}(\da_C/K)^{-1}\delta_K,
\end{equation}
where $K = \Fq(C)$ and $(\omega)$  is the divisor of a nonzero rational differential form $\omega \in \Omega^1_K$ on $C$.
There is the Plancherel formula $ \langle f, g\rangle =  \langle \Fo(f), \Fo(g)\rangle$ where  $f \in {\cal D}(\da_C)$,  $g \in {\cal D}'(\da_C)$ and $\langle -,-\rangle$ is the canonical pairing
between dual spaces.
When we apply this formula to  the characteristic functions  $ \delta_{\A_1(D)}$ and $\delta_K$ the result easily
yields  Riemann-Roch theorem together with  Serre duality for
divisors on  $C$ (see for example \cite{P6}).

Trying  to extend the harmonic analysis to the higher local fields and adelic groups we meet the following obstacle. The $n$-dimensional local fields
and consequently  the adelic groups are not locally compact topological groups
for $n>1$ in any reasonable sense whereas by a  theorem  of Weil   the existence of  Haar measure (in the usual
sense) on a topological group implies its local compactness. Unfortunately, the well-known extensions of this measure theory to the infinite-dimensional
spaces or groups (such as  the Wiener measure) do not help in our circumstances. Thus, we have to develop a measure theory and harmonic analysis on $n$-dimensional local fields and adelic groups {\em ab ovo}.

The idea for dealing with this problem  came to me in the  1990s. In dimension one,   local fields and adelic groups are equipped with a natural filtration provided by  fractional ideals $\wp^n,~n \in \Z$,  which correspond to the standard valuations. For example, this filtration on the field $\Fq ((t))$  is given by the powers of  $t$. If  $P \supset Q$ are two elements of such a  filtration on  a group $V$,  then the Bruhat space ${\cal  D}(V)$ is canonically isomorphic to the double inductive limit of the (finite-dimensional) spaces  ${\cal  F}(P/Q)$ of all functions on the finite groups $P/Q$.  The other function spaces listed above can be represented in the same way if we use all possible combinations of projective or inductive limits.

 In dimension two,   local fields $K$ such as $K_{P, C}$ again have  a  filtration by fractional ideals, which are  powers of $\wp$. But now, the quotient $P/Q = \wp^m/\wp^n, n > m$ will be isomorphic  to a direct sum of finitely many copies of the residue field $\bar K = \Fq ((u))$. Thus this group is locally compact  and the functional space    ${\cal  D}(P/Q)$ is well defined. To define the function spaces on $K$ one can try to repeat the procedure which we know for the 1-dimensional fields. To do that, we need to define the maps (direct or inverse images) between the spaces  ${\cal  D}(P/Q), {\cal  D}(P/R),{\cal  D}(Q/R)$  for $ P \supset Q \supset R$. The selection rules mentioned above restrict the opportunities for this construction. This enables us to introduce the following six types of spaces of functions (or distributions) on $V$:

$$
\begin{array}{lcccl}
{\cal D}_{P_0}(V) &=& \mbox{lim} & \mbox{lim} & {\cal D}(P/Q) \otimes \mu(P_0/Q) \mbox{,} \\
        &&  \stackrel{\longleftarrow}{j^*} & \stackrel{\longleftarrow}{i_*} & \\
{\cal D'}_{P_0}(V) &=& \mbox{lim} & \mbox{lim} & {\cal D}'(P/Q)  \otimes \mu(P_0/Q)^{-1} \mbox{,} \\
        &&  \stackrel{\longrightarrow}{j_*} & \stackrel{\longrightarrow}{i^*} & \\
{\cal E}(V) &=& \mbox{lim} & \mbox{lim} & {\cal E}(P/Q) \mbox{,} \\
&& \stackrel{\longleftarrow}{j^*} &
\stackrel{\longrightarrow}{i^*} & \\
{\cal E'}(V) &=& \mbox{lim} & \mbox{lim} & {\cal E}'(P/Q) \mbox{,} \\
&&  \stackrel{\longrightarrow}{j_*} &
\stackrel{\longleftarrow}{i_*} & \\
\tilde{{\cal E}}(V) &=& \mbox{lim} & \mbox{lim} & \tilde{{\cal E}}(P/Q) \mbox{,} \\
&&  \stackrel{\longrightarrow}{i^*}    &
\stackrel{\longleftarrow}{j^*} & \\
\tilde{{\cal E}}'(V) &=& \mbox{lim} & \mbox{lim} & \tilde{{\cal E}}'(P/Q) \mbox{,} \\
& & \stackrel{\longleftarrow}{i_*}
 &  \stackrel{\longrightarrow}{j_*}\\
&&&& \\
\end{array}
$$
where $P \supset Q \supset R$ are some elements of the filtration in
$V$ (with locally compact quotients), $P_0$ is a fixed subgroup from the
filtration and $j: Q/R \rightarrow P/R$, $i: P/R \rightarrow P/Q$
are the canonical maps.

This definition works for a general class of groups $V$ including the adelic groups such as $\da_X$, which has a
filtration by the subspaces $\da_{12}(D)$ where $D$ runs through the Cartier divisors
on $X$.

Thus, developing of harmonic analysis may  start with
the case of dimension zero (finite-dimensional vector spaces over a
finite field representing a scheme of dimension zero, such as $\Spec (\Fq)$, or finite abelian
groups) and then be extended by induction to the higher dimensions.

 An important contribution was made in 2001 by Michael  Kapranov \cite{K2} who  suggested using
a trick from the construction of the Sato Grassmanian in the theory of
integrable systems (known as a construction of  semi-infinite monomials)\footnote{A construction of this kind for the local fields is also contained in \cite{K}.}. The idea consists of  using
the spaces $\mu(P_0/Q)$ of measures instead of $\mu(P/Q)$ in the above definition of
the spaces ${\cal D}_{P_0}(V)$ and ${\cal D'}_{P_0}(V)$: without it one cannot define   the functional spaces
 for {\em all} adelic groups in the two-dimensional case and,
in particular, for the whole adelic space $\da_X$.

In 2005 Denis Osipov has introduced the  notion of a  $C_n$ structure in the category of filtered vector
spaces~\cite{O4}. With this notion at hand,   harmonic analysis can be developed in a very general setting, for all objects of the category
$C_2$. The crucial point is that the $C_n$-structure exists for the adelic spaces of {\em any} $n$-dimensional
noetherian scheme. The principal advantage
of this approach is that one can perform all the constructions {\em
simultaneously} in the  local and global cases. The category $C_1$
contains (as a full subcategory) the category of linearly locally
compact vector spaces (introduced and thoroughly studied by S.
Lefschetz \cite{L}) and there one can use the classical harmonic analysis.

When we go to  general arithmetic schemes over $\Spec (\Z)$,  fields like  $\C((t_1))\dots((t_n))$ appear and we need to extend the basic category $C_n$. In dimension one, this means that   connected Lie groups must also be considered. It is possible to define categories  of  filtered abelian groups  $C^{ar}_n,~(n = 0, 1, 2$), which contain all types of  groups which  arise from arbitrary  schemes of dimension 0, 1 and 2 (in particular  from algebraic surfaces over $\Fq$ and arithmetic surfaces).
Harmonic analysis can be developed  for these categories if  we introduce  function spaces which are close to  that of classical functional analysis,
such as  Schwartz space ${\cal S}(\R)$  of smooth  functions on  $\R$, which are rapidly  decreasing together  with all their derivatives.
 Recall that in the case of dimension one  we had to consider, in addition to   the genuine local fields such as $\Q_p$,  the fields $\R$ and $\C$. In the next dimension, we have to add to the two-dimensional local fields  such as $ \Fq ((u))((t))$ or  $\Q_p((t))$ the fields $\R((t))$ and $\C((t))$.
They will occupy the entire   row in  the table above. This  theory has been  developed in  papers \cite{OP1, OP2}.

Just as in the case of dimension one, we
define direct and inverse images in the categories of groups,
which take into account  all the components of the adelic complex, the Fourier transform
$\Fo$ which  preserves the spaces ${\cal D}$ and
${\cal D'}$ but interchanges the spaces
${\cal E}$ and ${\cal E'}$. We also introduce  the characteristic
functions $\delta_W$ of subgroups $W$  and then
prove a generalization of the Poisson
formula. It is important that for a certain class of
groups $V$ (but not for $\da_X$ itself) there
exists a nonzero  invariant measure, defined up to multiplication by a
constant, which is an element of ${\cal D'}(V)$.
Another important tool of the theory are the base change theorems for the inverse and direct images. They are function-theoretic
counterparts of the classical base change theorems in the categories of coherent sheaves.

The applications of the theory includes  an analytic expression for
the intersection number of two divisors based on an adelic approach to the
intersection theory \cite{P3} and  an analytic proof of the
(easy part of) Riemann-Roch theorem for divisors on $X$.

This theory is the harmonic analysis on the additive groups of the local fields and adelic rings (including their archimedean cousins).
 In the classical case of dimension one, the analysis can be
developed on arbitrary varieties (defined either over $K$, or over $\da$).
This has already been done by Bruhat in the local case \cite{Br}. For
arbitrary varieties defined over a two-dimensional local field $K$,
this kind of analysis  was carried out  by D. Gaitsgory and D.~A. Kazhdan in
\cite{GK1} for the purposes of representation theory of reductive
groups over the field $K$. This was preceded by a construction \cite{K3} of  harmonic analysis on
homogenous spaces such as $G(K)/G(O'_K)$ (introduced in \cite{P4}).
We note that the construction of harmonic analysis (over $K$ and $\da$) is
 a nontrivial problem  even in  the case $G = {\mathbb G}_m$. This will be  the
topic of our   discussion in the following sections.

\section{Discrete Adelic Groups on Two-dimensional Schemes}

The harmonic analysis discussed  above can be viewed as a representation theory of the simplest algebraic group over local or adelic rings,
namely, of  the additive group. In general, 1-dimensional local fields and adelic rings lead to a vastly developed representation theory of reductive groups
over these fields and rings. The simplest case of this theory is still the case of  an abelian group, namely  $\GL(1)$.   Let $K$ be a local field of dimension 1.
Then $\GL(1, K) = K^*$, the multiplicative group of $K$, and the irreducible representations are the abelian characters, i.e. continuous homomorphisms
$\chi : K^* \rightarrow \C^*$. For arithmetic applications one requires the morphisms to  $\C^*$, not to the unitary group ${\mathbb U}(1) \subset \C^*$.

The  1-dimensional local field $K$ contains a discrete valuation subring ${\cal O}$ with a maximal ideal $\wp$. Then the local group  $K^*$ has the following structure
$$
K^* = \{ t^n, n \in \Z \}\times {\cal O}^* = \{ t^n, n \in \Z \}\times {\bar K}^*  \times  \{1 +   \wp\} ,
$$
where $t$ is a generator of the ideal $\wp$, $\bar K = \Fq$ and the  group  $\{1 +   \wp\}$ is the  projective limit of its finite quotients  $ \{1 +   \wp\}/ \{1 +   \wp^n\}$. Thus, our group $K^*$  is a product of the maximal compact subgroup ${\cal O}^*$ and a  discrete group $\cong \Z$.
When $K$ is the local field $K_P$ attached to a point $P$ of an algebraic curve $C$ defined over a finite field $\Fq$,
 let us set  $\Gamma_P: = K_P^*/{\cal O}_P^*$.   In  the adelic case,  we set
$$
\Gamma_C: = \A_C^*/ \prod_P  {\cal O}_P^* = \bigoplus_P  K^*_P/{\cal O}_P^* = \bigoplus_P \Z.
$$
This group is the group of divisors on  $C$.

 We now introduce the  groups dual to these discrete groups viewing  them as algebraic groups defined over $\C$:
$$
\begin{array}{cccccccccccc}
{\mathbb T}_P &=& \Hom (\Gamma_P, \C^*), && {\mathbb T}_S &=& \prod_{P \in  S} {\mathbb T}_P, && {\mathbb T}_C& = &\lim &{\mathbb T}_S,\\
&&&&&&&&&&\stackrel{\longleftarrow}{S}&
\end{array}
$$
where $S$ runs through  all finite subsets in $C$. Let us consider the divisor $D_S$ with normal crossings  on  ${\mathbb T}_S$ that consists of the points in the product ${\mathbb T}_S$  for which  at least one component is the identity point in some ${\mathbb T}_P$. Let $\C_+[{\mathbb T}_S]$ be the space of  rational functions on ${\mathbb T}_S$ that are regular outside  $D_S$ and may have poles of  first order on $D_S$.
 The space  $\C_+[{\mathbb T}_C]$ can be defined as an inductive limit with respect to the obvious inclusions.

We would like to show that  harmonic analysis on the adelic space $\da_C$ can be reformulated in terms of complex analysis on the dual groups. We need one more torus  ${\mathbb T}_0 \cong \C^*$,  which corresponds  by the duality to  the image of the degree map
$$
\deg: \Gamma_C \rightarrow \Z\quad \mbox{with}~ \deg (D) = \sum_P n_P \deg (P)~\mbox{for a divisor}~ D = \sum_P n_P P.
$$

Denote by $j: {\mathbb T}_0  \rightarrow {\mathbb T}_C$ the natural embedding. Then the following diagram
\begin{equation}
\begin{CD}
{\cal D}(\da_C)^{{\cal O}^*}  = :{\cal D}_+(\Gamma_C)     @> {\cal L} >>   \C_+[{\mathbb T}_C]  @> j^* >>   {\cal F}_+[{\mathbb T}_0]\\
@ V \Fo VV     @ VVV    @ V i^*VV\\
{\cal D}(\da_C)^{{\cal O}^*}  = :{\cal D}_+(\Gamma_C)     @> {\cal L} >>   \C_+[{\mathbb T}_C]  @> j^* >>    {\cal F}_+[{\mathbb T}_0]
\end{CD}
\end{equation}
commutes. Here, the map $\Fo$  is induced by the Fourier transform on the adelic group $\da_C$, the map $i : {\mathbb T}_0 \rightarrow {\mathbb T}_0$
sends  $z \in {\mathbb T}_0$  to $q^{-1}z^{-1}$ and the space ${\cal F}_+[{\mathbb T}_0]$ consists of  the functions that are regular outside the
points $z = 1$ and $ z = q^{-1}$ and may have poles of the first order at these points.
We denoted here by  ${\cal L}$ a duality map, a version of the Fourier transform in this situation (completely different however from the Fourier map $\Fo$).
 If $g \in G$ and $z \in {\mathbb T}_G = \Hom (G, \C^*)$ for some group $G$ then $({\cal L}f)(z) = \sum_g f(g)z(g)$.

The next important fact is a reformulation of the Poisson formula on the group $\da_C$\footnote{For the sake of simplicity, we assume that $\Pic^0(C)(\Fq) = (0)$, that is $\Ker(\deg) = \Div_l(C)$.}. It can be shown that for any function $f \in {\cal D}(\da_C)^{{\cal O}^*}$
$$
\sum_{\gamma \in K}f(\gamma) = \res_{(0)}(\omega) +  \res_{(1)}(\omega),
$$
$$
\sum_{\gamma \in K}(\Fo f)(\gamma) = - \res_{(q^{-1})}(\omega) -  \res_{(\infty)}(\omega),
$$
where $\omega = j^*{\cal L}fdz/z$ is the differential form on the compactification of the torus  ${\mathbb T}_0$ and the points we have chosen for the residues are $z = 0, z = q^{-1}, z =1$ and $z = \infty$. Since the poles of the form $\omega$ are contained in this set,  we deduce that the Poisson formula on the curve $C$ (with an appropriate choice of  Haar measure on $\da_C$) is equivalent to the residue formula (4) for the form $\omega$ on the compactification
of the torus ${\mathbb T}_0$ (the general case see in \cite{P10}).

Our main goal now is to understand what correspond to these constructions in the case of dimension two\footnote{We consider here the case of an algebraic surface. The main definitions remain valid for the scheme part of an arithmetic surface.}.
Let us first consider the local situation, that is we fix a  flag $P \in C$ on $X$ and assume, for the sake of simplicity, that $P$ is a smooth point on $C$.
The local field $K_{P, C}$ has the  discrete valuation subring $\widehat{\cal O}_{P, C}$. It is mapped onto the local field
$k(C)_P$ on  $C$. This local field contains his own discrete valuation subring $\widehat{\cal O}_P$ and we denote its preimage in
$\widehat{\cal O}_{P, C}$ by $\widehat{\cal O}'_{P, C}$
We set
$$
\Gamma_{P,C}:=K_{P, C}^*/\widehat{\cal O}'^{\,*}_{P,C}
$$
where $\Gamma_{P, C}$ is a certain abelian group, which is (non-canonically) isomorphic to $\Z\oplus\Z$.
However,  there is a canonical exact sequence of abelian groups
\begin{equation}\label{equat-localex}
0\to\Z\to\Gamma_{P,C}\to\Z\to 0.
\end{equation}
The  map to $\Z$ in the sequence corresponds to the discrete valuation $\nu_C$ with respect to  $C$ and the subgroup $\Z$ corresponds to the discrete valuation $\nu_P$
on  $C$ at  $P$.
A choice of local coordinates $u, t$ in a neighborhood  of  $P$ such that locally $ C = \{ t = 0\}$
provides  a splitting of this exact sequence.
The group $\Gamma_{P,C}$ will then be isomorphic to the subgroup $\{ t^n u^m, n, m \in \Z \}$ in $K_{P, C}^*$.

The  group of coordinate transformations $u \mapsto u, t \mapsto u^kt, ~k \in \Z$  preserves  extension (\ref{equat-localex}).
Therefore,  this determines an embedding
\begin{equation}\label{aut}
\Z \to \Aut(\Gamma_{P,C}).
\end{equation}
which in fact is canonical.
\bigskip

 We are now going  to produce  a  global analogue of the local construction given above. For that purpose, consider  the subgroup $\widehat{\cal O}'^{\,*}$ of $\da^*_X$,  defined as the adelic product of the local groups  $\widehat{\cal O}'^{\,*}_{P,C}$ for all flags on an algebraic surface $X$.
Let us  consider the quotient
$$
\Gamma_X:=\da^*_X/\widehat{\cal O}'^{\,*}=:\prod_{(P, C)}{'\Gamma_{P, C}}.
$$
We have  a  natural surjective homomorphism $ \da^*_X\to\Gamma_X$ and all subgroups in  $\da^*_X$ such as $\da^*_{01}, \da^*_{12}, \dots,
\da^*_0$ have their images $\Gamma_{01}, \Gamma_{12}, \dots, \Gamma_{0} $ in   $ \Gamma_X$.

Then the structure of    $ \Gamma_X$ can be described by  an exact sequence
\begin{equation}\label{exact}
0\to \prod_{C}\Div(C)\longrightarrow\Gamma_X\stackrel{\pi}\longrightarrow
\bigoplus_{C}\prod_{P\in C}{'\Z}\to 0\,,
\end{equation}
where, as above, $\Div(C)$ denotes the group of  divisors on a curve $C \subset X$ and the restricted  product $\prod{'\Z}$ denotes the set of collections
of integers  with  components  whose absolute values are  bounded.
More precisely,
\begin{enumerate}
\item[(1)]
The subgroups $\prod_{C}\Div(C)$ and $\Gamma_{12}$ in $\Gamma_X$ coincide.
\item[(2)]
The restriction of the homomorphism $\pi$ to the subgroup $\Gamma_{02}\subset \Gamma_X$ is an isomorphism:
$$
\pi|_{\Gamma_{02}}:\Gamma_{02}\stackrel{\sim}\longrightarrow\bigoplus_{C}\prod_{P\in C}{'\Z}.
$$
\end{enumerate}

In other words, we see that there is a canonical splitting $\Gamma_X =  \Gamma_{12}\oplus\Gamma_{02}$ of  exact sequence (\ref{exact})
which is independent  of any possible choice of the coordinates.
The groups which we have constructed are  abelian. In our two-dimensional  case, the crucial point   is that they are provided with certain  canonical central extensions.

Let us  start once more  with the local situation, that is  we fix a  flag $P \in C$ on $X$.
Following \cite{AKdC}(see also \cite{KP1}) we have a   canonical central extension of groups
\begin{equation}\label{eq-AKC}
1\to k(C)^*_P \to {\tilde K}_{P, C}^* \to K_{P, C}^* \to 1\,.
\end{equation}
such that the corresponding commutator map
in the central extension is a skew form
$\langle\cdot,\cdot\rangle:K_{P, C}^*\times K_{P, C}^*\to k(C)^*_P$
given  by the tame symbol (without sign), that is by
\begin{equation}\label{ts}
\langle f,g \rangle=  f^{\nu_C(g)}g^{-\nu_C(f)}~({\rm mod}\,\wp)\in k(C)^*_P,
\end{equation}
where $\wp$ is the ideal  which defines the curve $C$.

There exists a canonical section of  extension~\eqref{eq-AKC} over the subgroup $\widehat{\cal O}'^*_{P, C}\subset K^*_{P, C}$. Denote
by $\tilde{\cal O}'^*_{P, C}$ the image of  $\widehat{\cal O}'^*_{P, C}$ in $\tilde{K}_{P, C}^*$ with respect to this section.
If we take the quotient of the extension ~\eqref{eq-AKC} by  the subgroup $\widehat{\cal O}^*_P$ of the center
$k(C)^*_P$ and then by the subgroup $\tilde{\cal O}'^*_{P, C}$   we obtain a  new central extension
\begin{equation}\label{eq-localext}
0\to \Z\to \tilde\Gamma_{P, C}\to\Gamma_{P, C}\to 0.
\end{equation}

 It is well known  that $H^2(\Z\oplus\Z,\Z)=\Z$ and the extension~\eqref{eq-localext} is a generator of this group. The commutator in this central extension defines a non-degenerate symplectic form $\langle -,-\rangle$ on  $\Gamma_{P, C}$ with values in $\Z$. Let us fix local parameters $u, t$ at $P$. Then
$\Gamma_{P, C}$ is isomorphic to the group of matrices
\begin{equation}\label{heis-matrix}
\left(
\begin{array}{ccc}
1&n&c\\
0&1&p\\
0&0&1
\end{array}
\right)
\end{equation}
with integer entries and $\langle n, p\rangle = np$. We denote this group by $\mbox{Heis}(3, \Z)$. Hence,  we arrive at  the following class of discrete nilpotent groups.

\smallskip

{\bf Definition 2}.
Let $H$, $H'$, and $\Co$ be abelian groups and let
$\langle-,-\rangle :H\times H'\to \Co$
be a biadditive pairing.  The set  $H \times H' \times \Co$ with the composition law $(n, p, c)(m, q, a) = (n + m, p + q, c + a + \langle n, q\rangle)$, where ãäå $n, m  \in H, p, q  \in H'$ and $c, a \in \Co$, is called  the {\em discrete Heisenberg group} $G$.

\smallskip

 One then constructs the Heisenberg group $G$ as a group of upper triangular unipotent matrices with $H$ and $H'$ on the second diagonal and $\Co$ in the right top corner. There is the obvious natural  central extension
$$
0\to \Co \to G\to H\oplus H'\to 0.
$$

In the global  case,  we have the Heisenberg group  $\tilde\Gamma_X$ with
$$
H=\Gamma_{12}=\prod_C\Div(C) = \prod_C\bigoplus_{P \in C}\Z,\quad
H'=\Gamma_{02}\cong \bigoplus_C\prod_{P \in C}{'\Z},
$$
$$
\Co=I_X:=\bigoplus_C  \bigoplus_{P \in C} \Z
$$
and the pairing $H\times H' \rightarrow \C$ is given by a component-wise multiplication.  We thus get a central extension
\begin{equation}\label{eq-extglobal}
0\to I_X\to \tilde\Gamma_X\to \Gamma_X\to 0
\end{equation}
and   for each flag $P \in C$  the restriction of  extension~\eqref{eq-extglobal} to $\Gamma_{P,C}$ coincides with  extension~\eqref{eq-localext}.
So, we obtain in this way    a global analogue  of the local construction, since we could  describe $\tilde\Gamma_X$ as an ``adelic'' product of the local groups $\tilde\Gamma_{P, C}$ in an appropriate  sense.

There is a natural surjective homomorphism $\varphi:I_X\to Z^2(X),\quad (n_{P, C})\mapsto\sum_{P}(\sum\limits_{C\ni P}n_{P,C})[P]$,
where $Z^2(X)$ denotes the group of zero-cycles on $X$. We set
$$
I_{02}:={\rm Ker}(\varphi),\quad I_{01}:=\bigoplus_C\Div_l(C)\subset \bigoplus_{C}\Div(C)=I_X.
$$
The Heisenberg group $\tilde\Gamma_X$ is closely related to the main arithmetic groups attached to  the surface $X$.
The quotient $I_X/(I_{01}+I_{02})$ is the second Chow group $CH^2(X)$ of $X$. Also, there are isomorphisms
$$
\Gamma_{01}/(\Gamma_0+\Gamma_1)\cong(\Gamma_{12}\cap(\Gamma_{01}+\Gamma_{02}))/\Gamma_1\cong
$$
$$
\cong(\Gamma_{02}\cap(\Gamma_{01}+\Gamma_{12}))/\Gamma_0\cong \Pic(X).
$$
Moreover, the pairing $\Gamma_{12}\times\Gamma_{02}\to I_X$ corresponds to the intersection pairing
$\Pic(X)\times\Pic(X)\to CH^2(X)$.

It is remarkable that the groups $K_{P, C}^*$ (and the  global adelic groups),
which are very far from being locally compact,   nevertheless have a non-trivial discrete quotient.
\section{Representations of Discrete Heisenberg Groups}
We have seen that in  the case of dimension two the first non-trivial nilpotent  groups  have occured. To define their duals one needs to develop an appropriate  representation theory for this class of groups.

 For  the discrete groups the classical theory of unitary representations on a Hilbert space  is not so well developed  since these groups are mostly  not of type I. By  Thoma's theorem,  a discrete group is of type I if and only if it has an abelian subgroup of finite index.

This implies a violation
of the main principles  of  representation theory on  Hilbert spaces: non-uniqueness of  the decomposition into  irreducible components;
too bad topology of the unitary dual space; non-existence of characters...  .
  V. S. Varadarajan wrote in 1989: ``A systematic developement of von Neumann's ideas led eventually (in the 1950s) to a deep understanding of the decomposition of unitary representations and to results  which implied more or less that a reasonable generalization of classical Fourier
analysis and representation theory could be expected only for the so-called {\em type I groups; i.e. groups all of whose factor representations  are of type I }"\cite{Var}.

 We can also say that the class of unitary representations is too restrictive for the arithmetic purposes.

On the other hand, there exists a theory of smooth representations for  $p$-adic algebraic groups. This theory is also valid for a more general class
of totally  disconnected locally compact groups. Discrete groups are a simple particular case of this class of groups and the general theory delivers a reasonable class of  representations, namely  representations on a vector space without any topology.
The new viewpoint consists in a systematic consideration of  purely algebraic  representations in place of   unitary representations on Hilbert spaces.

Following \cite{P11},  we  consider  now  this representation theory for the discrete Heisenberg groups $G = (H, H', \Co, \langle -, -\rangle)$ where all three groups are finitely generated.
We introduce the complex tori $\T_H = \Hom (H, \C^*),\T_{H'} = \Hom (H', \C^*)$  and  $\T_{\Co} = \Hom (\Co, \C^*)$, and set $\T_G = \T_H \times \T_{H'} \times \T_{\Co}$.
The group  $H$ is homomorphically mapped to $\T_{H'}$ according to the rule:
\begin{equation}\label{rule}
h \in H \mapsto  \{h' \mapsto \chi_{\Co}(\langle h, h'\rangle)\}.
\end{equation}

 Denote the kernel of this map by  $H_{\chi}$.
If  $\chi \in \T_{H'}$ then let $h(\chi)$ be the translate of the character  $\chi$  by the image of  $h$ in $\T_{H'}$. We have $h(\chi_{H'})(p) = \chi_{H'}(p)\chi_{\Co}(\langle h, p\rangle)$ for any $p \in H'$. For any $\chi \in \T_G, \chi = \chi_H \otimes \chi_{H'}\otimes \chi_C$, let $ G_{\chi} = H_{\chi}H'\Co$ be the subset in  $G$.
Then  $G_{\chi}$ is a normal subgroup in $G$, which depends only on $\chi_C$ and $\chi\vert G_{\chi}$ is a character of the group $G_{\chi}$ \cite{Pt}.

\smallskip

{\bf Definition 3}. Let  $V_{\chi}$ be the space of all complex-valued functions $f$ on $G$ which satisfy the following conditions:
\begin{enumerate}
\item $f(hg) = \chi(h)f(g)~ \mbox{for all}~ h \in  G_{\chi}$.
\item The support $Supp(f)$ is contained in the union of a finite number of  cosets of  $G_{\chi}$.
\end{enumerate}
Right  translations  define a representation $\pi_{\chi}$ of the group  $G$ on the space $V_{\chi}$.
 One can prove that these representations $\pi_{\chi}$ are irreducible in both possible senses: there are no nontrivial invariant subspaces, and the Schur lemma holds. Furthermore, these representations can be completely classified. Namely,
the representations  $V_{\chi}$ and  $V_{\chi'}$ are equivalent if and only if  three following conditions are satisfied:
\begin{enumerate}
\item $\chi_{\Co} = \chi'_{\Co}$.
\item There exists $h \in H$ such that  $\chi'_{H'} = h(\chi_{H'})$.
\item $\chi'_{H}(h) = \chi_{H}(h)$ for all $h \in H_{\chi}$ or equivalently there exists $t \in \T_{H/H_{\chi}} = \Hom(H/H_{\chi}, \C^*) \subset \T_H$
  such that $\chi'_{H} = t(\chi_{H})$.
\end{enumerate}

Here the torus $\T_{H/H_{\chi}}$ acts on the  ambient torus $\T_H$ by translations. The equivalence classes of representations
$V_{\chi}$ therefore correspond to  orbits of the groups  $\T_{H/H_{\chi}}\times H/H_{\chi}$  in  subsets  $\T_H \times \T_{H'} \times \{\chi_C\}$ of the torus $\T_G$ .

The group  $G$ is a semidirect product of the groups $H$ and $H'C$ and
the main tool for obtaining  the results stated above is
 the Mackey formalism \cite{M} which   describes the category of induced representations for semi-direct products of abelian locally compact groups. In the  classical theory, this is well-known for  unitary representations on Hilbert spaces. In our case, we can use the  version of this formalism developed in the theory of representation of $p$-adic reductive groups \cite{BZ, C, V}.

The restriction of functions from  the group  $G$ to the subgroup $H$  defines a bijection of $V_{\chi}$ with a certain space of functions on $H$. This space has an explicit  basis and we can now define the character of the representation $\pi_{\chi}$ as the matrix trace   of the representation operators $\pi_{\chi}(g)$ with respect to  this basis. It is easy to see that in many cases the corresponding infinite sum of  diagonal elements will  diverge.
 The simplest example is the group $\mbox{Heis}(3, \Z)$, see  (\ref{heis-matrix}).

 It is nevertheless possible  to define  the character if  we apply a well-known construction from the theory of loop groups \cite{PS}[ch. 14.1]. Namely, we have to add some "loop rotations "  to the group $G$ . In our context, this means that the group $G$ has to be extended to a semi-direct product  $\hat G = G \rtimes A$, where  $A \subset \Hom(H, H')$ is a non-trivial subgroup.

To construct the group $\hat G = G \rtimes A$, one needs to extend the automorphisms of the abelian groups $H\oplus H'$ to the automorphisms of the entire Heisenberg group. Note that the group $A$ acts on $H\oplus H'$ by  unipotent transformations. When  we fix an $r \in H$ and choose
 $k \in A$, the   expression
$$
k(m, p, c) = (m, p + k(m), c + 1/2\langle m - r, k(m)\rangle)~\mbox{ãäå}~m \in H, p \in H', c \in \Co
$$
 defines an automorphism of the group $G$ if the following conditions hold:
\begin{enumerate}
\item  $\langle m, k(m')\rangle~ =~ \langle m', k(m)\rangle~\mbox{for all}~m, m' \in H$
\item $\langle m - r, k(m)\rangle~ \in 2\Co~\mbox{for all}~m \in H.$
\end{enumerate}
In the case of the group  $\Gamma_{P, C}\cong \mbox{Heis}(3, \Z)$,  (see  (\ref{heis-matrix})),
we have $A = \Z, r = 1$ and $k \in A$ acts as $k(n, p, c) = (n, p + kn, c +\frac{1}{2}kn(n-1)), n, p, c \in \Z$.
Then the  extension is suggested by the existence of the group of coordinate transformations on the surface $X$ (see (\ref{aut})) at the point $P$:
$$
t \mapsto u^k t, ~u \mapsto u, ~k \in \Z.
$$
According to the analogy between algebraic and arithmetic surfaces we discussed above, these coordinate transformations in the  two-dimensional local field $K_{P, C} = \Fq ((u))((t))$ indeed correspond to the loop rotations
$$
t \mapsto \alpha t, ~\alpha \mapsto \alpha, ~\alpha \in \C
$$
in the field $\C ((t)$ on an arithmetic surface.

\smallskip

When $k(H_{\chi}) \subset \Ker (\chi_{H'})$ the representation of  $G$ on   $V_{\chi}$ can be extended to a representation $\hat\pi_{\chi}$ of the extended group $\hat G$ on  the same space.  Let
$$
(\T_{\Co} \times A)_+ : = \{\chi \in \T_{\Co}, k \in A : \mid \chi_{\Co}(\langle n, k(n)\rangle)\mid < 1~\mbox{for all}~ n \in H/H_{\chi}, n \ne 0 \}
$$
be a relation in $\T_{\Co} \times A$,  let $A(\chi)$  be the projection of the set $(\T_{\Co} \times A)_+ \cap (\{\chi\} \times A)$
to $A$ and let ${\hat G}(\chi) = G \times A(\chi) \subset {\hat G}$.

We can now solve the existence problem for the characters.
 The trace $\Tr~ \hat \pi_{\chi}(g)$ exists for all $g  \in {\hat G}(\chi)$
and we have
$$
\Tr~ \hat\pi_{\chi}(g) = \chi_{H}(m)\chi_{H'}(p)\chi_{\Co}(c)\cdot \sum_{n \in H/H_{\chi}} \chi_{H'}(k(n))\chi_{\Co}(\langle n, p\rangle +
1/2\langle n - r, k(n)\rangle).
$$
for  $g = (m, p, c, k),~k \in A(\chi),~m \in  H_{\chi} $. The trace is zero if  $m$ does not belong to $H_{\chi}$.

The trace is well-defined,  but  does not determine a function on the set of equivalence classes of  representations. To overcome this difficulty,
we have to consider  representations of the extended group ${\hat G}$.

Let $\T_A = \Hom(A, \C^*)$ and  $ \T_{\hat G} = \T_G \times \T_A$. If $\hat\chi = (\chi, \chi_A) \in \T_{\hat G}$, then we set
$$
\hat\pi_{\hat\chi} = \hat\pi_{\chi}\otimes \chi_A.
$$
We therefore have   $\Tr~ \hat\pi_{\hat\chi} = \Tr~ \hat\pi_{\chi}\cdot \chi_A$.
For a given $g \in \hat G(\chi)$,  the trace   $\Tr~ \hat\pi_{\hat\chi}(g)$ can be considered  as a function on the domain  $T' = \T_{H} \times \T_{H'}\times \T_{\Co}(k) \times \T_A$ in the torus  $\T_{\hat G}$,  where $\T_{\Co}(k)$  is the projection of the set $(\T_{\Co} \times A)_+ \cap (\T_{\Co} \times \{k\})$  to the torus $\T_{\Co}$.

Let us define an  action of the group  $\T_{H/H_{\chi}} \times H$ on the set  $\T_{H} \times \T_{H'}\times \{\chi_{\Co}\} \times \T_A \subset T'$ by the formula
\begin{equation}\label{action}
 (t, h) (\chi_{H}, \chi_{H'}, \chi_{\Co}, \chi_A) = (t(\chi_{H}), h(\chi_{H'}), \chi_{\Co}, \chi^{'}_{A}),
 \end{equation}
  where
 $$
  \chi^{'}_{A}(k) = \chi_A(k)\chi_{H'}(k(h))\chi_{\Co}(1/2\langle h - r, k(h)\rangle),~ k \in A.
 $$

We define the  space   ${\cal M}_G(k), k \in A$ as the quotient of the domain  $T'$ by this action.
The quotient-space is a complex-analytic manifold, in fact a fibration over a domain in  $\T_{\Co}$.
For a given $g = (m, p, c, k)\in {\hat G}(\chi)$ the trace  $Tr~ \hat\pi_{\hat\chi}(g)$ is invariant,  under a simple additional condition, under the action (\ref{action}) and  defines  a    holomorphic function $F_g =  F_g(\hat\chi)$ on ${\cal M}_G(k)$.  We now obtain  the main property that the characters must enjoy:

\bigskip

 Let  $\hat\chi ,\hat\chi^{'} \in \T_{\hat G}$. The representations    $\hat\pi_{\hat\chi}$ and  $\hat\pi_{\hat\chi^{'}}$  are equivalent if and only if  ${\hat G}(\chi) = {\hat G}(\chi')$ and  $F_{g}(\hat\chi) = F_{g}(\hat\chi')$ for all $g \in {\hat G}(\chi)$.

\bigskip

Thus we see that the space   ${\cal M}_G(k)$ is actually a moduli space for a  certain class of  representations of  ${\hat G}$.

Let us consider the simplest example, that  of the group $\mbox{Heis}(3, \Z)$.
Let $A = \Z = \Hom(H, H'),~r = 1,~ \hat G = G \rtimes \Z$ and $\chi_{\Co}(c) = \lambda^c,~\chi_{\Co} \in \T_{\Co}(k > 0)$ where $\T_{\Co}(k > 0)  = \{0 < |\lambda| < 1\}$. Then $\T_{H'}/\mbox{Im} H =: E_{\lambda}$  is an elliptic curve,  where $z \in \T_{H'} = \C^* ,~ \mbox{Im} H = \{ \lambda^{\Z}\}$. We have
a degree map
$$
\Pic (E_{\lambda}) = H^1(E_{\lambda}, {\cal O}^*) = H^1(H, {\cal O}^*(\T_{H'})) \rightarrow \Hom (H, H') = A,
$$
and
$$
\Pic (E_{\lambda}) = \{ \varphi(n, z) = a^{-n}z^{-kn}\lambda^{-1/2kn(n-1)}: a \in \C^*,~k \in A = \Z \}.
$$
Let $L$ be the line bundle which   corresponds to a 1-cocycle    $\varphi$. Then
$$
H^0(E_{\lambda}, L) = \{f(z), z \in \T_{H'}: f(\lambda^nz) = \varphi(n, z)f(z) \}.
$$
 The theta-series
$$
\vartheta_{p,k,a}(z, \lambda): = z^p\sum_{n \in \Z}a^nz^{kn}\lambda^{np + 1/2kn(n-1)}
$$
(which  are the Poincar\'e series with  respect to   $\varphi$)  converge  for all $z \in \C^*,~0 <  |\lambda| < 1, ~k >0$, and form
 a basis of the space  $H^0(E_{\lambda}, L)$ for $0 \le p < k$. Finally,
\begin{equation}\label{trace}
\Tr~ \hat\pi_{\hat\chi}(0, p, c, k) = \lambda^ct^k\vartheta_{p,k,1}(z, \lambda), ~(z, \lambda) \in {\cal A}_G(k),~t \in \T_A .
\end{equation}
In this case,   the theta-series lifted to  ${\overline{\cal A}}_G(k) = \C \times \{\mbox{upper  halfplane} \}$  are Jacobi modular forms  (up to some powers of $\lambda$ and $z$) with  respect to the standard  action of a finite index subgroup of the group  $(\Z\oplus \Z) \rtimes\mbox{SL}(2, \Z)$). This  last statement  is completely parallel to a well-known property of characters for  representations of affine Kac-Moody algebras \cite{KP2, PS}.

 In the more general situation in which   $H$ and  $H'$ are torsion-free groups and  $C = \Z$, $\vert \chi_C(c)\vert \neq 1$ for $c \neq 0$, the map $k: H \rightarrow H'$ is a monomorphism with finite cokernel,  $A = \Z k$ and the form $\langle -,k(-)\rangle$ is positive-definite, we have two dual abelian varieties  $E = \T_{H'}/\mbox{Im} H$ and  $E' = \T_{H}/\mbox{Im} H'$ with the Poincar\'e bundle ${\cal P}$ over $E\times E'$.
  The morphism $k$ defines an isogeny $\varphi_k : E \rightarrow E'$ and the sheaf $L$ is defined as $(Id\times \varphi_k)^*{\cal P}$.
By Mumford's theory  \cite{Mum}, there exists a finite Heisenberg group $\widetilde{\Ker} (\varphi_k)$, which  is a central extension of the group $\Ker (\varphi_k)$.
Then for all   $g = (m, p, c, k)\in {\hat G}(\chi)$  the values of the characters $\Tr~ \hat\pi_{\hat\chi}(g)\chi_{\Co}^{-1}(c)\chi_A^{-1}(k)$ are theta-functions for  the bundle $L$.

 If $\hat\chi = 1\otimes \chi_{H'}\otimes \chi_{\Co}\otimes 1$,  then the functions $\Tr~ \hat\pi_{\hat\chi}(0, p, 0, k)$  for  $p \in H'\mbox{mod}~k(H)$ form a basis of  the space    $H^0(E, L)$. This basis is  a standard Mumford basis for the action of the  Heisenberg group   $\widetilde{\Ker} (\varphi_k) = (H'/H, \T_{H'/H}, \C^*)$    on the space   $H^0(E, L)$.

 In addition, certain orthogonality relations are satisfied by the characters \cite{P11}.

The boundary of the domain  $\T_{\Co}(k)$ can contain those  characters  $\chi_0 \in \T_{\Co}$ for which  $H_{\chi_0}$ has a finite index in $H$. These characters correspond to the roots of unity in   $\C^*$, so that the  representations $\pi_{\chi_0}$ are finite-dimensional.
Let $V = H\otimes \R$ and $Q$ be the extension of the pairing $\langle n, k(n) \rangle ,~n \in H$ to the space $V$. Also, let $\chi_{\Co}(c) = \lambda^c$ and let us choose  a boundary point $\chi_0$.
The classical limit formulas for  theta-functions imply the following behavior of the trace near the  $\chi_0$ (we assume that $\chi_H = 1$ and $\chi_H' = 1$):
\begin{equation}\label{limit}
\Tr~\hat\pi_{\hat\chi}(g) \sim \Tr~\hat\pi_{\hat\chi_0}(g) \cdot [H : H_{\chi_0}]^{-1}(\Det_V Q)^{-1/2}(\frac{\sqrt \pi}{2})^{\rk H}\log {\vert \lambda \vert}^{-\frac{1}{2}\rk H}\quad \mbox{when}~\chi_{\Co} \to \chi_0.
\end{equation}
The trace of the representation $\hat\pi_{\hat\chi_0}$ can be computed in terms of a Gauss sum.

Thus, we see that, in our situation, the change in  the class
of representations  will cause  the moduli spaces of induced representations to be   complex-analytic manifolds.
Characters do exist and are the  modular forms. It seems that this more general holomorphic dual
space is more adequate for this class of groups than the standard unitary dual which goes  back to the Pontrjagin duality for abelian groups.

\section{Problems and Perspectives}
 We collect here several problems related to the issues we have discussed in the talk.

\smallskip

{\bf 1.} {\bf Harmonic analysis for  local fields and adelic groups of arbitrary dimension ${\bf n}$}.

The basic category for this study has to be the category   $C_n$ \cite{O4} and its version that  includes  fields of the archimedean type \cite{OP2}. When one tries to extend the  measure theory and harmonic analysis to   $n$-dimensional local fields and adelic groups for $n > 2$ the following problem arises. The selection rules become too severe to go further in a straightforward way.  This obstacle appears  already for  $n$-dimensional local fields with $n = 3$. We can define the spaces analogous to  ${\cal D}(V)$ or ${\cal D}'(V)$ only under some strong restrictions on the groups $V$ (= objects in $C_n$). Note  that   spaces such as ${\cal E}(V)$ can be easily defined for any $n$ and arbitrary group $V$.

\smallskip

{\bf 2}. {\bf The Tate-Iwasawa method for two-dimensional schemes.}

\smallskip

J. Tate \cite{T1}  and independently K. Iwasawa \cite{I}  reformulated the classical problem of analytic continuation
for zeta- and L- functions for the   fields of algebraic numbers and the  fields of algebraic functions in one variable over a finite field.
They  introduced a new type of   $L$-functions:
$$
L(s,\chi, f)=\int_{\da^*} f(g)\chi(g)|g|^s d^*g
$$ where $d^*g$ is a Haar measure on $\Bbb A^*$, the function $f$ belongs to the Bruhat-Schwartz space of
functions on $\Bbb A_X$ and $\chi$ is an abelian character of the group $\Bbb A^*$ associated to a character
$$\chi: \Gal(K^{\mbox{ab}}/K) \rightarrow \C^*
$$
of the Galois group  by the reciprocity map $\Bbb A^* \rightarrow \Gal(K^{\mbox{ab}}/K)$ . They also  proved the
analytic continuation of  $L(s,\chi, f)$ to the entire $s$-plane and the functional equation
$$
L(s, \chi, f) = L(1-s,\chi^{-1}, \Fo(f))
$$
by means of  the Fourier transform $\Fo$ and the Poisson formula for  functions  on  $\Bbb A_X$ (\ref{p1}), (\ref{p2}).

For a special
choice of $f$ and $\chi = 1$ we obtain  the zeta-function
$$
\zeta_X(s) = \prod_{x \in X} (1 - (\#k(x))^{-s})^{-1},
$$
of any  scheme $X$ of dimension one
(to which we have to add, if necessary,  the archimedean factors).
Here $x$ runs through the closed points of  $X$.
The product converges for $\text{Re}(s) > \dim X$.

There exists a general Hasse-Weil conjecture \cite{Ha, W1}
which asserts that these zeta- (and more general $L$-) functions can be meromorphically extended to the
entire  $s$-plane and satisfy  the functional equation (for regular proper schemes $X$ of dimension $n$)
of the type $\zeta_X(n - s) = \{\mbox{elementary factors}\}~\zeta_X(s)$.

This conjecture has been  completely proved for  algebraic varieties defined over a finite field $\Fq$.  For this goal the powerful machinery of
the \'etale cohomology has been  developed by A. Grothendieck.  For  schemes over $\Spec (\Z)$, the general results are
known only in dimension one,  thanks to  the Hecke's theorem. Later  this was included into the  Tate-Iwasawa approach. At the same time,  this approach works for algebraic curves defined
over $\Fq$. For the higher dimensions over $\Spec (\Z)$, there are only  scattered results; however these include the
proof of the Hasse-Weil conjecture for elliptic curves over $\Q$ \cite{Wi, BCDT}.

 For a long time the author has advocated the following

{\bf Problem.} {\em Extend Tate--Iwasawa's analytic method to higher dimensions} (see in particular \cite{P6}).

The higher adeles were  introduced precisely for this purpose. We hope that  harmonic analysis and  representation theory of  adelic groups on two-dimensional schemes may help to solve this problem.

\smallskip

{\bf 3.} {\bf  Behavior of zeta- and L-functions in the critical strip.}

\smallskip

The critical strip for the ordinary Riemann's zeta-function is  $0 \le \Re (s) \le 1$ and this zeta-function (with an archimedean factor) has there exactly two poles, both  of  first order. For the two-dimensional case, the critical strip is wider, namely $0 \le \Re (s) \le 2$. Take as  $X$ a model over $\Spec (\Z)$ of an elliptic curve $E$ defined over $\Q$. The Birch and Swinnerton-Dyer conjecture \cite{BSD, T2} states  that
\begin{equation}\label{bsd}
  \zeta_X (s) \underset{s\to 1}{\sim} \frac{\# E(\Q)_{tor}^2}{c~ \Omega~ \Det_{E(\Q)}\langle -, -\rangle~ \#\Sh } (s - 1)^{-r - 2},
\end{equation}
where $E(\Q)$ is the finitely generated Mordell-Weil group of rational points on $E$, $r$ is its rank, $\langle -, -\rangle$ is the height pairing, $\Omega$ is the real period of the curve, $\Sh$  is the  Shafarevich-Tate group and $c$ is a product of certain local invariants.

Many years ago several people, including the author, have independently observed that this limit behavior is very similar to the limit behavior of a theta-function attached to a lattice. Namely, let $V/\R$ be a finite dimensional euclidean vector space of dimension $n$.
Denote by $\langle -, -\rangle$ the scalar product on  $V$. Let $\Gamma $ be a finitely generated abelian group such that $\Gamma \otimes \R = V$ and
let $\Gamma' = \Gamma /\Gamma_{tor}$ be the corresponding  lattice (= a discrete co-compact subgroup) in $V$.
Then the theta-function $\theta_\Gamma(t)$ is defined as
$$ \theta_\Gamma(t)\ :=\
      \sum_{\gamma\in \Gamma} e^{-\pi t \langle \gamma, \gamma \rangle} =  \#\Gamma_{tor}\cdot \theta_{\Gamma'}(t)
$$
and satisfies the functional equation
$$ \theta_{\Gamma'}(t)\ =\
     t^{-\frac{n}{2}} \, \mbox{Vol}(\Gamma')^{-1} \theta_{\Gamma'^\perp}(t^{-1}) \ ,
$$
where $\Gamma'^\perp\subset V$ is the dual lattice and the volume of the fundamental domain for $\Gamma'$ is
$\mbox{Vol}(\Gamma') = \det(\langle e_i, e_j\rangle)^{1/2}$ with $\{e_i\}$ a basis of
 the free $\Z$-module $\Gamma'$.

In particular, we get
$$
 \theta_\Gamma(t)\ \underset{t\to 0}{\sim}\
         \#\Gamma_{tor} \mbox{Vol}(\Gamma)^{-1} t^{-\frac{n}{2}}\ .
$$
If we apply this asymptotic  formula to the group $\Gamma \oplus \Gamma$
then we  get
\begin{equation}\label{limit1}
 \theta_{\Gamma \oplus \Gamma}(t) \ \underset{t\to 0}{\sim}\
      \frac{\# \Gamma_{tor}^2}{\Det_{\Gamma}\langle -, - \rangle} t^{- n},
\end{equation}
which looks rather similar to the conjecture (\ref{bsd}) if we take as   $\Gamma$  the group $ E(\Q) \oplus \Z \oplus \Z$. D.  Zagier has devoted to this relation a note \cite{Z}  with many interesting remarks and observations. In particular, he discussed the question of interpreting     such factors   as $\Omega$ and $\# \Sh$  which are not visible in the theta-formula (\ref{limit1}).

In order to clarify the situation,  let us  look at the corresponding behavior of the zeta-function of an algebraic surface $X$ defined over $\Fq$.
The analogy between geometric surfaces over $\Fq$ and arithmetic surfaces such as this model $X$ of $E$ suggests that this may be a useful
move.

The value of the zeta function at $s=1$ is given by the conjecture of Artin and Tate \cite{T2,  Mil}.
 We assume that  $X$ is  a smooth proper irreducible  surface.
 Denote by $\rho=\rk~ \mbox{NS}(X)$ the rank of the Neron-Severi-group of $X$ and
 let $\{D_i\}$ with $D_i\in \mbox{NS}(X)$ $i=1,\ldots,\rho$ be a basis of
 $\in \mbox{NS}(X)\otimes \Q$. Denote by $D_i \cdot D_j$ their intersection index.
 Let $\mbox{Br}(X)=H^2(X_{et},{\cal O}_X)$ be the Brauer group of $X$.
 Then the group $\mbox{Br}(X)$ is conjectured to be finite  and the following relation holds:
 $$
     \zeta_X (s) \underset{s\to 1}{\sim}  (-1)^{\rho-1}\,q^{\chi({\cal O}_X)}\,
       \frac{\# \Pic(X)_{tor}^2}{\# H^0(X, {\cal O}^*_X)^2 \,\# \mbox{Br}(X)\,\det((D_i\cdot D_j))}
         \left( 1\,-\,q^{1-s} \right)^{-\rho}\ .
  $$
Within the framework of the analogy between geometry and arithmetic \cite{P9}, the group $\mbox{NS}(X)$ corresponds to the group  $E(\Q) \oplus \Z \oplus \Z$,
the intersection index corresponds to the height pairing, the period $\Omega$  corresponds to  $q^{\chi({\cal O}_X)}$ and the Brauer group to the Shafarevich-Tate  group $\Sh$.

Since $ \left( 1\,-\,q^{1-s} \right)^{-\rho} \underset{s\to 1}{\sim} (s - 1)^{-\rho}(\log q)^{-\rho}$, we  again guess that  certain
theta-functions related to the lattice $\mbox{NS}(X)$  may have this kind of the limit behavior.
An immediate objection to this suggestion is that the intersection pairing is not positive-definite. This  can be resolved if we consider the Siegel theta-functions attached to indefinite quadratic forms.

The case of surfaces makes it clear that this  question is highly non-trivial. Zeta-functions of algebraic varieties over $\Fq$ are very simple analytic functions. Indeed, according to  Grothendieck's  theory, they are equal to $F(q^{-s})$ where $F(t)$ is a rational function of a variable $t$.  The theta-functions involved  are certainly transcendental  functions, which  cannot be  simplified in this way by  substitution.  Thus the problem we arrive at  is to understand how theta-functions can appear in this setting in a natural way,  and how to relate them to zeta-functions. We conjecture  that the theta-functions which occur into the traces of representations of the adelic groups constructed above could be such theta-functions. Their  behavior in the limit (\ref{limit}) has the structure we have just described.

It is worth  mention  another problem, the so called  $S$-duality conjecture, which is quite close to what have been  discussed here.
The problem came from the quantum field theory \cite{VW} but has purely algebraic formulation for an algebraic  surface $X$ over a finite field $\Fq$ (see a discussion in \cite{K1}). Let $M_{r, n}$ be a moduli space of semi-stable vector bundles $E$ on $X$ with given rank $r$, trivial determinant and the second Chern class $c_2(E) = n$.
Then the formal series
$$
\sum_n  \#M_{r, n}(\Fq) q^{-ns}
$$
is expected to have  under mild conditions on $X$ a modular behavior with respect to  a congruence subgroup of the group $\mbox{SL}(2, \Z)$. It is remarkable that the transcendental functions  appear once more in relation to a surface defined over a finite field.

\smallskip

{\bf 4. Representations of  discrete nilpotent groups.}

\smallskip

i) The representations  $\pi_{\chi}$ and  $\hat\pi_{\hat\chi}$ of the discrete Heisenberg groups are  particular examples of the irreducible  representations of these groups.  Thus, the  problem of classification of all irreducible  representations arises. Of course, one needs to impose certain  conditions in order to get a reasonable answer. In the theory of unitary representations  for discrete nilpotent finitely generated groups $G$ on a Hilbert space such a condition was found in \cite{B}. One says that a representation $\pi$ of  $G$ on a space $V$ has the finite multiplicity property if there exists  a subgroup $H \subset G$ which  preserves a line $l$ in $V$ and such that the character of $H$ defined by the action of $H$ on $l$ occurs in $\pi\vert_H$ as a  discrete direct summand with finite multiplicity.
Then the class of irreducible representations with this property coincides with the class of irreducible monomial (= induced by an one-dimensional character) representations of $G$.

It is highly  desirable to define in our algebraic situation a class of ``basic'' induced  representations which will play the  role  that the  Verma modules or  representations with  highest weight do for  the representations of reductive Lie groups (or algebraic groups). This is closely related to a problem of classification of (say, left) maximal ideals in the group ring of  $G$.

ii)  The moduli spaces ${\cal M}_G(k)$ defined above are  orbit spaces for group actions. This construction looks very similar to the Kirillov's orbit method for connected real (or complex) nilpotent Lie groups $G$ (or nilpotent algebraic groups over $\Q_p$) \cite{Kir} where the unitary dual is the space ${\mathfrak g}^*/G$ of co-adjoint orbits in the dual ${\mathfrak g}^*$ of the Lie algebra ${\mathfrak g}$ of $G$.   Attempts to extend Kirillov's method to   finitely generated nilpotent groups were  made in \cite{H, Koh} (see also \cite{BS}). It seems that there is a general functorial definition of spaces such as  ${\cal M}_G(k)$ for  arbitrary nilpotent discrete groups which will replace the spaces  ${\mathfrak g}^*/G$ in this situation, just as  the torus ${\mathbb T}_{\hat G}$ may be an analogue of  the space  ${\mathfrak g}^*$. The Kirillov's character formula may also exist in this situation.

iii) When one  tries to apply  the representation theory developed in section 5 to  the nilpotent groups which  arise from the algebraic surfaces $X$ (section 4), one immediately observes  that:

1) the groups like $\tilde\Gamma_X$ are not finitely generated;

2) the groups like (\Pic (X), \Pic (X), \mbox{CH} (X)) are  equipped with the indefinite form $\langle -, -\rangle$.

 Certainly, the representation theory cannot be automatically extended to the case of  infinitely generated  groups. In our case, the ``big" group
$\tilde\Gamma_X$ is the  adelic product of  simplest Heisenberg groups $\tilde\Gamma_{P, C}$ and  consequently  is an inductive limit of  finite products of these local groups. We can easily extend  all the representation-theoretic  constructions to the case of    $\tilde\Gamma_X$ if we apply the technique from  the theory of adelic products of  reductive algebraic groups over 1-dimensional local fields. The role of the compact subgroups is now played by  co-finite products of the local Heisenberg groups.

The problem  2)  can also be  solved. A solution is based on using   the Siegel  theta-functions for indefinite
quadratic forms that are well suited for this situation.

iv) An  important problem is to develop an analysis on discrete Heisenberg  groups $G$, in particular, to define appropriate  function spaces on $G$ , the analogue of the map ${\cal L}$ (see (8) in section 3) and to obtain a Plancherel-type theorem  which  relates  the function spaces on $G$  and spaces of holomorphic (or meromorphic) functions on  ${\cal M}_G(k)$.

v) There exists a  general question of the decomposition into the irreducible components  of representations of  discrete nilpotent groups. It is known that the regular representation  (on the $L^2$-space on $G$) of a discrete group $G$ may have  very different decompositions into irreducible components (see a first example of this kind in \cite{M}). On the other hand,  in our situation there is a rather concrete problem:  how does one  decompose the natural  fundamental representation  of the group $\tilde\Gamma_X$ (and locally  of the groups $\tilde\Gamma_{P, C}$)  on  the spaces  ${\cal D}_{\da_{12}}(\da_X)^{{\cal O}'^*}$  or ${\cal D}_{\da_{12}}'(\da_X)^{{\cal O}'^*}$(respectively in ${\cal D}_{{\cal O}_{P, C}}(K_{P, C})^{{\cal O}'^*_{P, C}}$ or  ${\cal D}_{{\cal O}_{P, C}}'(K_{P, C})^{{\cal O}'^*_{P, C}}$) on a surface $X$ ?

vi)  Our  theory deals with  the discrete ``part" of the adelic group $\da_X^* = \GL (1, \da_X)$. D. Gaitsgory and D. Kazhdan  have extended  the traditional theory of representations for  reductive $p$-adic groups (parabolic induction, Jacquet functor, cuspidal representations) to the case of groups $\GL (n, K)$ where $K$ is a two-dimensional local field (and of more general reductive groups)\cite{GK1, GK2, GK3}.  An  important and certainly very hard problem is to merge these two theories, at least for the group $\GL (2, \da_X)$.

vii) For the schemes of dimension two, we constructed  discrete Heisenberg groups,  which are  nilpotent groups of  class 2. It is possible to associate certain  discrete adelic groups  to schemes of arbitrary dimension $n$ and that are the  nilpotent groups of  class $n$.

\bigskip

In this text, we mainly gave   a review of  certain recent advances in the  higher adelic  theory.
During  the last thirty years, this  theory was developed in many different directions. We finish with  a short list of these developements\footnote{This list of references does not pretend to be complete.}:
\begin{itemize}
\item residues and symbols \cite{P1, P2, FP, Y,  BM1, BM2, BM3, Kh1, Kh2, O2,  Sop,  OZ}
\item  class field theory  for higher dimensions: the author,  K. Kato and his school, S. V. Vostokov and his school, see surveys \cite{FP, FV, Inv, R}
\item adelic resolutions for sheaves, intersection theory, Chern classes, Lefschetz formula for coherent sheaves \cite{P3, Y,  HY1, HY2, GP, G}
\item  algebraic groups over local fields, buildings, Hecke algebras \cite{P4, P8, K3, GK1, GK2, GK3,  BK}
\item restricted adelic complexes and the Krichever correspondence \cite{P7, O1, O3, KOZh1, KOZh2}
\item relations with non-commutative algebra \cite{P5,  Zh}.
\end{itemize}

\end{document}